\documentclass[12pt]{article}

\usepackage[utf8]{inputenc}

\usepackage{amsfonts}
\usepackage{amssymb}
\usepackage{amsmath}
\usepackage{amsthm}
\usepackage[T2A]{fontenc}

\usepackage[russian]{babel}
\usepackage{amsthm}

\def \le {\leqslant}
\def \ge {\geqslant}

\theoremstyle{plain}

\begin{document}

\begin{Large}
\centerline{\bf \"Uber eine Ungleichung von Schmidt und Summerer}
\centerline{ \bf f\"ur diophantische Exponenten }
\centerline{\bf von Linearenformen in drei Variablen }
\end{Large}
\vskip+1cm
\begin{large}
\centerline{von Nikolay Moshchevitin\footnote{ Die Untersuchung ist  von der Beihilfe der russischen Regierung 
11.
G34.31.0053 und  RFBR  No.12-01-00681-a unterstützt.
 } (Moskau)}
 \end{large}
\vskip+2cm

{\bf \S 1. Diophantische Exponenten.}

Sei  $n\in\mathbb{Z}_+$ und $\Theta=(\theta_1,...,\theta_n)\in \mathbb{R}^n$.
Das $n$-Tupel  $\Theta$ hei\ss t  eigentlich, wenn $1,
\theta_1,...,\theta_n$  linear unabh\"angig sind  \"uber dem  K\"orper der rationalen Zahlen.
 Ist ein $n$-Tupel $\Theta$ gegeben,
so setze man f\"ur $t\ge 1$
$$
\psi_\Theta (t) = \min_{x_1,x_2,x_3 \in \mathbb{Z}:\,\, 1\le \max_j |x_j|\le t}||\theta_1x_1+\theta_2x_2 + \theta_3 x_3||.
$$
Die {\it diophantischen Exponenten} $\omega = \omega(\Theta)$ und $\hat{\omega} = \hat{\omega}(\Theta)$ sind definiert durch:
$$
\omega = \omega(\Theta) =\sup \{\gamma\in \mathbb{R}:\,\,
\liminf_{t\to+\infty} t^\gamma \psi_\Theta (t) < +\infty\},
$$
$$
\hat{\omega} = \hat{\omega}(\Theta)
=
\sup \{\gamma\in \mathbb{R}:\,\,
\limsup_{t\to+\infty} t^\gamma \psi_\Theta (t) < +\infty\}.
$$
 Bekanntlich ist
$$
\psi_\Theta (t) \le t^{-n},
$$
sodass 
\begin{equation}\label{u}
\omega(\Theta )\ge \hat{\omega}(\Theta)  \ge n
.
\end{equation}

Im Falle $n=1$ haben wir $ \hat{\omega}(\Theta) = 1 $  f\"ur alle $\Theta\in \mathbb{R}\setminus \mathbb{Q}$.

Im Falle $n=2$ hat 
V. Jarn\'{\i}k [1] eine Verbesserung der trivialen Ungleichung  (\ref{u}) erhalten.

{\bf Satz 1.}  (V. Jarn\'{\i}k [1]) 
{\it
Ist
$\Theta=(\theta_1,\theta_2)$
eigentlich, so ist
\begin{equation}\label{1q}
\omega(\Theta) \ge \hat{\omega}(\Theta) (\hat{\omega}(\Theta)+1).
\end{equation}
}
M. Laurent 
[3]  zeigt da\ss \, die Ungleichung (\ref{1q}) optimal ist.

Im Falle $n\ge 3$ zeigt V. Jarn\'{\i}k [2] eine weitere Ungleichung, die aber nicht optimal ist. Es wurde eine Verbesserung von Moshchevitin [4,5] vorgenommen, die sich aber auch als nicht optimal erwiesen hat. 

W.M. Schmidt und L. Summerer  [7,8]
haben eine neue und leistungsf\"ahige Methode entwickelt, mit ber sie die folgende Ungleichung bewiesen haben.

{\bf Satz 2.} (W.M. Schmidt und L. Summerer  [9])\,\,{\it Ist
$\Theta=(\theta_1,\theta_2,\theta_3)$
ein 
eigentliches System, so ist
\begin{equation}\label{ua}
\omega (\Theta ) \ge \hat{\omega}(\Theta)\cdot \frac{\sqrt{4\hat{\omega}(\Theta) - 3}-1}{2}.
\end{equation}
}

Satz 2 enthält nun die optimale Ungleichung.

Weiter möchten wir Satz 2 kurz beweisen, indem wir uns auf unsere \"Uberlegungen in [6] beziehen.

\vskip+0.3cm

{\bf  \S 2. Ein kurzer Beweis für Satz 2.}

 Es seien
$$
 {\bf m}_\nu = (m_{0,\nu},m_{1.\nu},m_{2,\nu}, m_{3,\nu})\in \mathbb{Z}^4$$ 
die {\it Bestapproximationsvektoren} f\"ur die lineare Form mit Koeffizienten $\theta_1,\theta_2,\theta_3$.
 Sei
$$
L_\nu = |m_{0,\nu}+m_{1.\nu}\theta_1+m_{2,\nu}\theta_2+ m_{3,\nu}\theta_3|,
\,\,\, M_\nu = \max_j |m_{j,\nu}|.
$$
Sei $\alpha <\hat{\omega}$, so hat man
\begin{equation}
L_{j-1}\le M_j^{-\alpha}
\label{1}
\end{equation}
wenn  $j$ hinreichend  gro\ss      \, ist.

Wir haben zwei F\"alle.

{\bf Fall 1.} Es gibt einen linearen Teilraum ${\cal L} \subset \mathbb{R}^4, {\rm dim} {\cal L} = 3$, und ${\bf m}_\nu \in {\cal L}$ wenn $\nu$
hinreichend gro\ss \, ist.
Im diesem Falle haben wir (\ref{1q}).  Das ist  besser  als (\ref{ua}).

{\bf Fall 2.} Es gibt unendlich viele  Paare $ \nu<k$  mit  den folgenden  Eigenschaften:

{\bf (i)} die Vektoren
$$
{\bf m}_{\nu-1},{\bf m}_\nu,{\bf m}_{\nu+1}
$$
sind linear unabh\"angig;

{\bf (ii)} die Vektoren
$${\bf
m}_{k-1},{\bf m}_k,{\bf m}_{k+1}
$$
sind linear unabh\"angig;

{\bf (iii)} 
es gibt  einen linearen Teilraum $\pi$, ${\rm dim }  \pi = 2$ und
 $$
{\bf m}_l\in \pi,\,\,\, \nu\le l\le k;\,\,\,\,\, {\bf m}_{\nu-1}
\not\in \pi,\,\,\, {\bf m}_{k+1} \not\in \pi;
$$

{\bf (iv)} die Vektoren
$$
{\bf m}_{\nu-1},{\bf m}_\nu,{\bf m}_{\nu+1},{\bf m}_{k+1}
$$
sind linearisch unabh\"angig.

Aus {\bf (iii)} folgt 
\begin{equation}\label{2}
L_{\nu}M_{\nu+1}\asymp L_{k-1}M_k.
\end{equation}Aus Eigenschaft {\bf (iv)} wird geschlossen auf
\begin{equation}\label{3}
1\ll L_{\nu-1}M_\nu M_{\nu+1}M_{k+1}.
\end{equation}
 Aus Ungleichungen 
 (\ref{2}) und (\ref{3}) ist herzuleiten
 $$
L_{\nu}M_{\nu+1}\asymp
L_{k-1}M_k \ll
L_{\nu-1}M_\nu M_{\nu+1}\cdot L_{k-1}M_kM_{k+1}.
$$
Sei
$$
b =
\frac{-1+\sqrt{4\alpha-3}}{2(\alpha-1)},\,\,\, a = 1-b.\,\,\,\,
a,b \in [0,1],\,\,\,\lambda = b(\alpha- 1) = \frac{\sqrt{4\alpha - 3}-1}{2}.
$$
Also entweder
\begin{equation}\label{a}
(L_{k-1}M_k)^a\ll L_{k-1}
M_kM_{k+1},
\end{equation}
oder
\begin{equation}\label{b}
(L_{\nu}M_{\nu+1})^b\ll L_{\nu-1}
M_\nu M_{\nu+1}.
\end{equation}
 Aus (\ref{a}) und (\ref{1}) folgt
$$
M_{k+1}\gg M_k^{b(\alpha -1)} = M_k ^\lambda
$$
und
\begin{equation}\label{p1}
L_k \ll M_k^{-\alpha\lambda}.
\end{equation}
Aus (\ref{b}) und (\ref{1}) hat man 
$$
L_\nu^b \ll M_{\nu+1}^{1-b}M_\nu^{1-\alpha}.
$$
Aber aus (\ref{1}) erh\"alt man
$$
M_{\nu+1}\le L_\nu^{-1/\alpha}
$$
und
$$
L_\nu^b \ll L_\nu ^{\frac{b-1}{\alpha}}M_\nu^{1-\alpha},
$$
sodass
\begin{equation}\label{p2}
L_\nu \ll 
M_\nu ^{\frac{\alpha(1-\alpha)}{b(\alpha -1)+1}} = M_\nu^{-\alpha\lambda}.
\end{equation}
Also  wir
 haben (\ref{p1}) oder (\ref{p2}). 

Daher gilt $\omega (\Theta ) \ge \lambda\alpha$,
und daraus folgt (\ref{ua}).

\vskip+0.3cm

Der Autor dankt  L. Summerer und Yu. Rabotkin für ihre Hilfe in der Übersetzung
 
\vskip+0.3cm
{\bf Literaturverzeichnis}

\noindent [1]\,
V.  Jarn\'{\i}k,\,\,
Une remarque sur les approximation diophantiennes
lin\'{e}aries,
Acta Scientarium Mathem. Szeged, 12 (pars B) (1949), 82 - 86.

\noindent [2]\,
В. Ярник,\,\,
К теории однородных линейных диофантовых приближений,
Чехословацкий математический журнал, т. 4 (79) (1954), 330 - 353.

\noindent [3]\, M. Laurent,\,\,
Exponents of Diophantine approximations in dimension two,
Canad. J. Math.  61, 1 (2009),165 - 189.

\noindent [4]\,
N.G. Moshchevitin,\,\ Contribution to Vojt\v{e}ch
Jarn\'{\i}k, Preprint available at arXiv:0912.2442v3 (2009).

\noindent [5]\,
N.G. Moshchevitin,\,\, Khintchine's singular Diophantine systems and their
applications, Russian Mathematical Surveys. 65:3(2010), 433 - 511.
 
\noindent [6]\,
N.G. Moshchevitin,\,\, Exponents for three-dimensional
simultaneous Diophantine approximations,
Czechoslovak Mathematical Journal, 62 (137) (2012), 127–137.

\noindent [7]\,
W.M. Schmidt, L. Summerer,\,\,
Parametric geometry of numbers and applications,  Acta Arithmetica, 140:1 (2009), 67-91.

\noindent [8]\,
W.M. Schmidt, L. Summerer,\,\, Diophantine approximation and parametric geometry of numbers,
Monatshefte f\"ur Mathematik, 169 (2013), 51 - 104.

\noindent [9]\,
W.M. Schmidt, L. Summerer,\,\,
Simultaneous approximation to three numbers,
Moscow J. Comb. Number Theory, 3:1 (2013), 84 - 107.

\vskip+1cm

Nikolay Moshchevitin

Fakult\"at f\"ur Mathematik und  Mechanik

Universit\"at Moskau 

Leninskie Gory

119991 Moskau

Ru\ss land

\vskip+0.2cm

e-mail: moshchevitin@gmail.com, moshchevitin@mach.math.msu.su

\end{document}